\documentclass[12pt]{amsart}
\usepackage{amssymb}
\usepackage{color}
\begin{document}

\def\dbl{[\hskip -1pt[}
\def\dbr{]\hskip -1pt]}
\newcommand{ \al}{\alpha}
\newcommand{\ab}[1]{\vert z\vert^{#1}}
\title
{
  Nonlinearizable CR automorphisms for polynomial models in
$\mathbb C^N$.
}
\author{martin kol\'a\v r and francine meylan}

\address{M. Kol\'a\v r: Department of Mathematics and Statistics, Masaryk University,
Kotlarska 2,  611 ~37
 Brno, Czech Republic}

\email{mkolar@math.muni.cz}
\thanks{The first author was supported by the GACR grant GA17-19437S}

\address{ F. Meylan: Department of Mathematics,
University of Fribourg, CH 1700 Perolles, Fribourg}

\email{francine.meylan@unifr.ch}

\begin{abstract}
We classify  polynomial models for real hypersurfaces 
in $\mathbb C^N$,  which admit nonlinearizable
infinitesimal CR automorphisms.  
As a consequence,
this provides an optimal 1-jet determination result in the general case.
Further we prove that  such automorphisms arise from one common source,
by pulling back  via  a holomorphic mapping a suitable symmetry of
a hyperquadric in some complex space. 

\end{abstract}

\def\Label#1{\label{#1}}
\def\1#1{\ov{#1}}
\def\2#1{\widetilde{#1}}
\def\6#1{\mathcal{#1}}
\def\4#1{\mathbb{#1}}
\def\3#1{\widehat{#1}}
\def\K{{\4K}}
\def\LL{{\4L}}
\def\H{{\4H}}
\def\C{{\4C}}
\def\R{{\4R}}
\def \MM{{\4M}}
\def\La{\Lambda}
\def\la{\lambda}
\def\Re{{\sf Re}\,}
\def\Im{{\sf Im}\,}

\numberwithin{equation}{section}
\def\s{s}
\def\k{\kappa}
\def\ov{\overline}
\def\span{\text{\rm span}}
\def\ad{\text{\rm ad }}
\def\tr{\text{\rm tr}}
\def\xo {{x_0}}
\def\Rk{\text{\rm Rk\,}}
\def\sg{\sigma}
\def \emxy{E_{(M,M')}(X,Y)}
\def \semxy{\scrE_{(M,M')}(X,Y)}
\def \jkxy {J^k(X,Y)}
\def \gkxy {G^k(X,Y)}
\def \exy {E(X,Y)}
\def \sexy{\sccenikrE(X,Y)}
\def \hn {holomorphically nondegenerate}
\def\hyp{hypersurface}
\def\prt#1{{\partial \over\partial #1}}
\def\det{{\text{\rm det}}}
\def\wob{{w\over B(z)}}
\def\co{\chi_1}
\def\po{p_0}
\def\fb {\bar f}
\def\gb {\bar g}
\def\Fb {\ov F}
\def\Gb {\ov G}
\def\Hb {\ov H}
\def\zb {\bar z}
\def\wb {\bar w}
\def \qb {\bar Q}
\def \t {\tau}
\def\z{\chi}
\def\w{\tau}
\def\Z{\zeta}
\def\phi{\varphi}
\def\eps{\varepsilon}

\def \T {\theta}
\def \Th {\Theta}
\def \L {\Lambda}
\def\b {\beta}
\def\a {\alpha}
\def\o {\omegacenik}
\def\l {\lambda}

\def \im{\text{\rm Im }}
\def \re{\text{\rm Re }}
\def \Char{\text{\rm Char }}
\def \supp{\text{\rm supp }}
\def \codim{\text{\rm codim }}
\def \Ht{\text{\rm ht }}
\def \Dt{\text{\rm dt }}
\def \hO{\widehat{\mathcal O}}
\def \cl{\text{\rm cl }}
\def \bR{\mathbb R}
\def \bS{\mathbb S}
\def \bK{\mathbb K}
\def \bD{\mathbb D}
\def \bC{\mathbb C}
\def \C{\mathbb C}
\def \N{\mathbb N}
\def \bL{\mathbb L}
\def \bZ{\mathbb Z}
\def \bN{\mathbb N}
\def \scrF{\mathcal F}
\def \scrK{\mathcal K}
\def \mc #1 {\mathcal {#1}}
\def \scrM{\mathcal M}
\def \cR{\mathcal R}
\def \scrJ{\mathcal J}
\def \scrA{\mathcal A}
\def \scrO{\mathcal O}
\def \scrV{\mathcal V}
\def \scrL{\mathcal L}
\def \scrE{\mathcal E}
\def \hol{\text{\rm hol}}
\def \aut{\text{\rm aut}}
\def \Aut{\text{\rm Aut}}
\def \J{\text{\rm Jac}}
\def\jet#1#2{J^{#1}_{#2}}
\def\gp#1{G^{#1}}
\def\gpo{\gp {2k_0}_0}
\def\emmp {\scrF(M,p;M',p')}
\def\rk{\text{\rm rk\,}}
\def\Orb{\text{\rm Orb\,}}
\def\Exp{\text{\rm Exp\,}}
\def\Span{\text{\rm span\,}}
\def\d{\partial}
\def\D{\3J}
\def\pr{{\rm pr}}

\def \CZZ {\C \dbl Z,\zeta \dbr}
\def \D{\text{\rm Der}\,}
\def \Rk{\text{\rm Rk}\,}
\def \CR{\text{\rm CR}}
\def \ima{\text{\rm im}\,}
\def \I {\mathcal I}

\def \M {\mathcal M}
\def\5#1{\mathfrak{#1}}

\newtheorem{Thm}{Theorem}[section]
\newtheorem{Cor}[Thm]{Corollary}
\newtheorem{Pro}[Thm]{Proposition}
\newtheorem{Lem}[Thm]{Lemma}
\newtheorem{Conj}[Thm]{Conjecture}

\theoremstyle{definition}\newtheorem{Def}[Thm]{Definition}

\theoremstyle{remark}
\newtheorem{Rem}[Thm]{Remark}
\newtheorem{Exa}[Thm]{Example}
\newtheorem{Exs}[Thm]{Examples}

\def\bl{\begin{Lem}}
\def\el{\end{Lem}}
\def\bp{\begin{Pro}}
\def\ep{\end{Pro}}
\def\bt{\begin{Thm}}
\def\et{\end{Thm}}
\def\bc{\begin{Cor}}
\def\ec{\end{Cor}}
\def\bd{\begin{Def}}
\def\ed{\end{Def}}
\def\br{\begin{Rem}}
\def\er{\end{Rem}}
\def\be{\begin{Exa}}
\def\ee{\end{Exa}}
\def\bpf{\begin{proof}}
\def\epf{\end{proof}}
\def\ben{\begin{enumerate}}
\def\een{\end{enumerate}}

\newcommand{\zz}{(z,\bar z)}

\maketitle

\section{Introduction}

 One of the fundamental questions in the theory of several complex variables, going back to the seminal work of H.~Poincar\'e (\cite{Po}),
is how to classify domains up to biholomorphic equivalence.
In the complex plane, the classical
Riemann mapping theorem
asserts 
that domains posses  only topological invariants.
 As was realized already by Poincar\'e, there is no analogous statement in higher dimension, 
and smooth boundaries of domains 
have infinitely many local biholomorphic invariants.

The main topic of this paper concerns  the classification  of   polynomial models of real hypersurfaces 
according to their infinitesimal CR automorphisms. 
The first  motivation for this study  comes from the fact that, as has been shown in \cite{KMZ},
the classical Chern-Moser theory can be extended to the case of singular Levi form - hypersurfaces with polynomial models of finite Catlin multitype. In particular, 
 it has been shown in \cite{KMZ} that  the kernel of the generalized Chern-Moser operator, which is 
  in one to one correspondence with the Lie algebra of infinitesimal CR automorphisms of the polynomial model of such  a hypersurface,
  gives a precise description of the derivatives needed to characterize an element of its stability group.
  
  In order to develop this approach further, towards a complete normal form construction and solution of the Poincar\'e local biholomorphic equivalence problem for such manifolds,
we need to classify the polynomial models according to the Lie algebra of infinitesimal CR automorphisms.

The second motivation comes from the study of possible sources and forms of  symmetries (infinitesimal CR automorphisms) of CR manifolds. 
Since linear  symmetries are relatively well understood,  
the main interest lies in understanding the possible existence and origin of nonlinearizable symmetries, provided by vector fields with vanishing linear part. 
This research  has 
long history, starting with the classical case of Levi nondegenerate hypersurfaces (\cite{CM},\cite{V}).
In this case, nonlinearizable automorphisms of the model hyperquadrics are
determined by their two jets. 

We consider this problem in the singular Levi form case. 
In this paper, we study systematically
nonlinear infinitesimal CR automorphisms of  {
polynomial models} in complex dimension {$N>3$,  the $\mathbb C^2$ and $\mathbb C^3$ cases beeing  completely understood} (see {\cite{Ko1a}}, \cite{KM2}).
The results provide a  description  of
hypersurfaces of finite Catlin multitype  in $\Bbb C^{N}$ 
{whose polynomial models admit such symmetries.}
 In combination with
the results of \cite{KMZ}, we prove a sharp 1-jet determination result for the holomorphic automorphism groups in general.

Moreover, we  identify the common source of such symmetries.
In the case of homogeneous polynomial models, they arise
via suitable holomorphic  mappings into   hyperquadrics in some complex spaces,  as  ``pull-back'' of  symmetries of the hyperquadrics.

 We first recall the sharp results of \cite{KMZ}
which describe explicitely the possible structure of the Lie algebra of infinitesimal CR automorphisms.

Let us consider
 a  holomorphically nondegenerate weighted homogeneous
model
of finite Catlin multitype in {$\Bbb C^{n+1},$} and denote
\begin{equation}\Label{mode0}
M_H:=\{\Im w = P_C(z,\bar z)\}, \quad (z,w)\in \C^n\times\C,
\end{equation}
where $P_C$ is a weighted homogeneous polynomial of degree one with respect to the multitype weights in the sense of Catlin
(for precise definitions, see Section 2).

It has been proved in \cite{KMZ}, that the   Lie algebra of infinitesimal CR automorphisms
$\5g=\aut(M_H,0)$ of $M_H$ admits the following weighted decomposition,
\begin{equation}\label{sid2}
\5g = \5g_{-1} \oplus \bigoplus_{j=1}^{s}\5g_{- \mu_j} \oplus \5g_{0}
\oplus 
 \5g_c \oplus  \5g_{nc}
 \oplus \5g_{1},
\end{equation}
where $\5g_{c}$ contains  vector fields commuting with
$W = \d_w$ and $\5g_{nc}$ contains vector fields not commuting with $W$, arising by  (possible) integration of nontransversal shifts. In both cases, 
the weights of such vector fields lie
in the interval $(0,1)$.
Note that apriori, $\5g_{nc}$ is only determined modulo $\5g_{c}$. However, it will be shown in the proof of Theorem \ref{th12} that there always exists a canonical 
representation for 
 $\5g_{nc}$. Note that $\5g_{c}=\{0\}$ in the case of a Levi nondegenerate hypersurface.
 Recall  that  the elements  of $\5g_{1}$  are the (possible) $2-$integrations of $W,$  and that  $\dim{\5g_{1}}=0$ or $1;$
 vector fields in $\5g_j$ with $j<0$ are regular and vector fields
in $\5g_0$ are linear (see\cite{KMZ}).

We introduce the following definition.
\bd
Let $R$  be a weighted homogeneous polynomial and $S$ be a holomorphic  weighted homogeneous polynomial. We 
say that $R$ admits an $S-$reproducing field, if there exists a holomorphic weighted homogeneous vector field  $Z$ such that
$$Z(R)=SR.$$
\ed
Note that for $S = 1$ we obtain the definition of a complex reproducing field, used in \cite{KMZ}
 \bt \label{th12}
Let $P_C(z,\bar z)$ be a weighted homogeneous polynomial of degree
1 with respect to the multitype weights,
such that the hypersurface
\begin{equation}\Label{mode0}
M_H:=\{\Im w = P_C(z,\bar z)\}, \quad (z,w)\in \C^n\times\C,
\end{equation} is holomorphically nondegenerate. Let
$\5g_{nc}$
in \eqref{sid2} satisfy \begin{equation}\dim \5g_{nc}>0.\end{equation}

Then   $M_H$ is biholomorphically
equivalent to

\begin{equation}\label{ex1}
\Im w =   \vert z_1 \vert^2 + S(z', \bar z') 
\end{equation}
or
{\begin{equation}
\Im w =  \Re ( z_1 \overline{ Q_1(z',\bar z')}) + S(z', \bar z'),
\label{ex2}
\end{equation}}
{where $Q_1$ is a holomorphic polynomial. }

Moreover, when $P_C$ is a  weighted homogeneous polynomial given by  \eqref{ex1}, then \eqref{ex1} admits a nontrivial $\5g_{nc}$ if and only if $S$ 
admits a reproducing field. 

Similarly, when $P_C$ is a homogeneous polynomial given by  \eqref{ex2} (i.e. 
the multitype weights are all equal), then  \eqref{ex2} admits a nontrivial $\5g_{nc}$ if and only if 
$Q_1$ and
 $S$  admit a common $Q_1-$ reproducing field $Y$, hence $Y(S) = Q_1S$ and  $Y(Q_1) = Q_1^2.$
\et
The  explicit description of the second case will be given  in Proposition \ref{third-nt};
in the  case of $\Bbb C^3$, we have   $S=0$ in the singular Levi form case. (See \cite{KM2}).
Note also  that the classical Levi nondegenerate case is covered by Theorem \ref{th12}.

We will show in Section 4 by an example (see Example \eqref{ex5}) that the last part of the claim does not hold in the case of unequal weights.

In order to describe hypersurfaces with nontrivial $\5g_c$ (which occur only in the singular Levi form case),  we introduce the following definition.

\bd
Let $Y$ be a weighted homogeneous vector field. 
A pair of finite sequences   of  vector valued holomorphic  weighted homogeneous polynomials  of dimension $s$ $\{U^{(1)}, \dots, U^{(l)}\}$ and  $\{V^{(1)}, \dots, V^{(l)}\}$ is called a symmetric pair of $Y_{s}-$chains if 
\begin{equation}
Y(U^{(l)})=0, \ \ Y(U^{(j)}) =A_j U^{{(j+1)}}, \ \ j=1, \dots, l-1,
\label{yuv}
\end{equation}
\begin{equation}
Y(V^{(l)})=0, \ \ Y(V^{(j)}) =B_j V^{{(j+1)}}, \ \ j=1, \dots, l-1,
\end{equation}
where $Aj$ and  $B_j$ are invertible $s \times   s$ matrices, which satisfy
\begin{equation} A_j= - ^t\overline {B_{l-j}}. \label{dcj}
\end{equation}
If the two sequences are identical,  
we say that $\{U^{(1)}, \dots, U^{(l)}\}$ is a symmetric $Y_s$ - chain.
\ed

The following result shows that in general the elements of $\5g_c$  arise from pairs of chains. 

\bt \label{pivo} Let 
 $M_H$  be  a holomorphically nondegenerate model  given by  \eqref{mode0} admitting a nontrivial 
 $Y\in \5g_c$. 
Then $P_C$ can be decomposed in the following way

\begin{equation}P_C=\sum_{j=1}^M T_j,\label{Barou190}
\end{equation} where each $T_j$ is given by 
\begin{equation}\label{Barou200}
 T_j =\Re (\sum_{k=1}^{N_j}<
  {U_j^{{(k)}}},
  { {{V_j^{{(N_j -k +1)}}}}}>), 
  \end{equation}
where  $\{{ {{U_j^{(1)}}, \dots, {U_j^{(N_j)}} }}\}$  and  $\{{ {{V_j^{(1)}}, \dots, {V_j^{(N_j)}} }}\}$ are  symmetric pairs of $Y_{s_j}-$ chains,  and $<\ ,\,>$ is the usual scalar product in $\Bbb C ^{s_j}.$

Conversely, if $Y$ and $P_C$ satisfy  \eqref{yuv} -- 
\eqref{Barou200},
then $Y \in \5g_c$.
\et


\bd  If $P_C$ satisfies \eqref{yuv} --
\eqref{Barou200}, the associated hypersurface  $M_H$ is called a chain hypersurface.
\ed

The description of the remaining component $\5g_1$ is a consequence of Theorem 4.7 in \cite{KMZ}
(see section 2 for the notations).
\bd We say that $P_C$ given by  \eqref{mode0} is balanced if it can be written as
\begin{equation}
 P_C(z, \bar z) = \sum_{|\al|_{\Lambda}= |\bar \al|_{\Lambda}
 = 1
 } A_{\al, \bar \al} z^{\al} \bar z^{\bar \al},
 \label{bala}
\end{equation}
for some nonzero n-tuple of real numbers  $\Lambda = (\lambda_1,\dots, \lambda_n)$,
where
$$\vert \alpha\vert_{\La} := \sum_{j=1}^n\la_j \al_j .$$
The associated hypersurface $M_H$ is called a balanced hypersurface.

\ed
Note that $P_C$ is balanced if and only if there exists a complex reproducing field $Y$ in the terminology of \cite{KMZ}, i.e.,
$ Y(P_C) = P_C$. Indeed, it  can be shown that such a $Y$ is of the form 
$$Y=\sum_{j=1}^n\lambda_j z_j \partial_{z_j} .$$ (See Lemma 4.6 in  \cite{KMZ}).

As a consequence of  Theorem 1.1 in \cite{KMZ}, we obtain the following result.
\bt\label{cap} The component
$\5g_1$ satisfies $\dim \5g_{1}>0 $  if and only if in suitable multitype coordinates $M_H$ is a balanced hypersurface.
\et
The following theorem gives    the number of  derivatives needed  to uniquely  determine the elements  of the stability group  of a class of smooth hypersurfaces in terms of their model hypersurfaces.

 \bt\label{raf}
 Let $M$ be a smooth hypersurface and $p \in M$ be a point of finite Catlin multitype with holomorphically
 nondegenerate model,  where $P_C$ is a  homogeneous polynomial.  If its model at $p$ is
 neither a balanced hypersurface nor a chain hypersurface,
 then its  automorphisms are determined by  the  1-jets at $p$.
 \et
Our last result is the following theorem

\bt\label{emb}
Let $M_H$ be a holomorphically nondegenerate hypersurface
given by \eqref{mode0}, where $P_C$ is a homogeneous polynomial and $Y$ be a vector field of strictly
positive weight.  Then $Y \in \aut(M_H, 0)$ 
{if and only if} there
exists an integer $K \ge n+1$ and a holomorphic mapping $f$ from a
neighbourhood of the origin in $\mathbb C^{n+1}$ into $\mathbb C^K$
which maps  $M_H$ into a Levi nondegenerate hyperquadric
$Q\subseteq \mathbb C^K$ such that the following holds:
\begin{enumerate}
 \item   $Y$ is $f$-related with  a $1-$integration of a  nontransversal shift of $Q$ if $Y \in\5g_{nc},$ 
 \item $Y$ is $f$-related with a $2-$integration of a  transversal shift of $Q$ if $Y \in\5g_{1},$ 

\item   $Y$ is $f$-related with a rotation of $Q$  if $Y \in \5g_c.$
\end{enumerate}
\et

The already mentioned example \eqref{ex5} suggests that in the case of unequal weights, (1) in Theorem \ref{emb} fails, although we do not prove this. 

Let us remark that mappings of CR manifolds into hyperquadrics have been studied intensively in recent years (see e.g.
\cite{BEH}, \cite{EHZ}).
Here we ask in addition that the mapping be compatible with a symmetry of the hyperquadric.

The paper is organized as follows. Section 2 contains the necessary definitions used in the rest of the paper.
Section 3 deals with the $\5g_{nc}$ component of the Lie algebra
$\5g$. Section 4 deals with the $\5g_c$  component while Section 5 contains the proofs of Theorem \ref{th12}, Theorem \ref{pivo}, Theorem \ref{cap}, Theorem \ref{raf} and 
Theorem \ref{emb}.

\section{Preliminaries}

In this section we recall the notion of Catlin multitype and  some definitions needed in the sequel.

Let  $M \subseteq \mathbb C^{n+1}$ be a smooth hypersurface,
and $p \in M $ be a  point of {\em finite type} $m$ in the
sense of Kohn and Bloom-Graham (\cite{BG}).
We will consider
local holomorphic coordinates $(z,w)$ vanishing at $p$,
where $z =(z_1, z_2, ..., z_n).$ 
 The hyperplane $\{ \Im w=0 \}$ is assumed to be tangent to
$M$ at $p$, hence  $M$  is described near $p$ as the graph of a uniquely
determined real valued function
\begin{equation} \Im w = \psi(z_1,\dots, z_n,  \bar z_1,\dots,\bar z_n,  \Re w), \ d\psi(p) \neq 0.
\label{vp1}
\end{equation}
Using a result of \cite{BG}, we may assume
that
\begin{equation}\label{fifi}
\psi(z_1,\dots, z_n,  \bar z_1,\dots,\bar z_n,  \Re w)=P_m(z, \bar z) +o(|\Re w|+|z|^m),
\end{equation}
where $P_m(z, \bar z)$ is a nonzero homogeneous polynomial of degree $m$ with {no pluriharmonic} terms.

The definition of multitype involves  {\it rational}  weights associated  to the variables
$w, z_1, \dots z_n$ in the following way.

  The
variables $w$, $\Re w,$ and $\Im w$ are given weight one, {reflecting} our choice
of  variables given by \eqref{vp1}.
The complex tangential variables $(z_1, \dots, z_n)$  are treated according to
the following definitions (for more details, see \cite{Ko1}, \cite{KMZ}).

\bd
A weight is an $n$-tuple of nonnegative
 rational numbers $\La = (\la_1, ...,
\la_n)$, where $0 \leq\la_j\leq \frac12$, and $\la_j \ge
\la_{j+1}$.
\ed
Let $\La = (\la_1, ...,
\la_n)$ be a weight, and   $\al=(\al_1,\dots, \al_n), \ $
 $\ \beta=(\beta_1, \dots,\beta_n) $   be  multiindices.
The weighted degree $\kappa$ of a monomial $q(z, \bar z,\Re w)=c_{\al
\beta l}z^{\al}\bar z^\beta {(\Re w)}^{l} , \ l \in \mathbb N,$ is defined as
$ \kappa:=
l +  \sum_{i=1}^n (\al_i + \beta_i ) \la_i.$

A polynomial $Q(z, \bar z, \Re w)$  is $\La-$homogeneous
of weighted degree $\kappa$ if it is a sum of
 monomials of weighted degree $\kappa$.

For a weight  $\La$,  the weighted length of a multiindex $\al = (\al_1, \dots, \al_n)$ is
defined by

$\vert \alpha\vert_{\La} := \la_1 \al_1 + \dots + \la_n \al_n.$

Similarly, if $\al = (\al_1, \dots, \al_n)$ and  $\hat \al =
(\hat \al_1, \dots, \hat \al_n)$ are two multiindices, the weighted
length of the  pair $(\al, \hat \al)$ is
$\vert (\alpha,\hat \al) \vert_{\La} := \la_1 (\al_1 +\hat \al_1) \dots +
\la_n (\al_n + \hat \al_n).$

The weighted order $\kappa$  of a differential operator is defined in a similar way.

\bd{ A weight $\La$ will be called distinguished for $M$ if there exist
local holomorphic coordinates $(z,w)$ in which the defining equation of $M$ takes form
\begin{equation} \Im w = P\zz + o_{\La}(1),
\label{1}
\end{equation}
where $P\zz$ is a nonzero $\La$ - homogeneous polynomial of
weighted degree $1$ without pluriharmonic terms, and $o_{\La}(1)$
denotes a smooth function whose derivatives of weighted order less than or equal to
one vanish.}
\ed
The  following definition is due to D. Catlin (\cite{C}).
\bd \Label{2.6} (\cite{C})
Let  $\Lambda_M = (\mu_1, \dots, \mu_n)$  be the infimum of all possible
distinguished weights  $\Lambda$ with respect to the lexicographic order.
The multitype of $M$ at $p$ is defined to be the $n$-tuple $(m_1,
m_2, \dots, m_n),$ where
$m_j = \begin{cases}   \frac1{\mu_j} \ \  {\text{ if}} \ \  \mu_j \neq 0\\
  \infty \ \ {\text{ if}} \ \   \mu_j = 0.
\end{cases} $

Furthermore, if none of the $m_j$ is
infinity, we say that $M$ is of {\it finite multitype at $p$}.

Coordinates corresponding to the multitype weight $\Lambda_M$, in
which the local description of $M$ has form (\ref{1}), with $P_C:=P$
being  $\Lambda_M$-homogeneous, are called
{\it multitype coordinates}.
\ed
 Note that if  $n=1,$ $M$ is of finite type at $p$ if and only if
$M$  is of finite multitype at $p.$ In this case, the type of $M$ at $p$  is equal to the multitype of $M$ at $p.$

From now on,  we assume that $p \in M$   is a point of {\it finite
multitype}, with $M$ defined locally by  

\begin{equation} \Im w = 
P_C\zz + o_{\La_M}(1).
\label{vp12}
\end{equation}

\bd (\cite{Ko1}){ Let $M$ be given by  (\ref{vp12}).
We define a
 model hypersurface $M_H$ associated
to $M$ at $p$ by
\begin{equation} M_H = \{(z,w) \in \mathbb C^{n+1}\ | \
 \Im w = P_C \zz \}. \label{2}\end{equation}}
\ed

Note that multitype coordinates $(z,w)$  are not unique. Nevertheless
it  is shown in \cite{Ko1} that all models are biholomorphically equivalent (in fact  by a polynomial transformation).

We recall that  $\Aut (M,p)$   is the stability group of $M,$ that is,  the set of those germs at $p$
of biholomorphisms mapping $M$ into itself and fixing $p,$
and  that  $\aut (M,p)$  is the set of  germs of holomorphic vector fields in
$\mathbb C^{n+1}$ whose real part 
is tangent to $M$.
If $M$ admits a holomorphic vector field $X$ in $\aut (M,p)$ such that $Im X $ is also tangent
(i.e.$X$ is complex tangent),
then $\aut (M,p)$ is of infinite
 dimension (\cite{S}).
\bd(\cite{S}){A real-analytic hypersurface $M \subset \mathbb C^{n+1}$  is  {\it holomorphically nondegenerate
at $p \in M$} if there is no germ at $p$ of a holomorphic vector field $X$ tangent to $M.$}
 \ed

 \bd{We say that the vector field $$Y = \sum_{j=1}^n F_j(z,w)
\partial_{z_j} + G(z,w)  \partial_{w} $$
 has homogeneous weight $\mu \  (\ge -1)$ if $F_j$ is  a weighted homogeneous polynomial
 of weighted degree $\mu+ \mu_j,$ and $G$ is a
 homogeneous polynomial of weighted degree $\mu+1.$}
\ed

 \bd\Label{types}
Let $X \in \aut(M_P,p)$ be a rigid  weighted homogeneous vector
field, that is, a vector field   whose coefficients do not depend on the $w$ variable.  $X$ is called
\begin{enumerate}
\item a {\it shift} if the weighted degree of $X$ is less than zero;
\item a {\it rotation }
if the weighted degree of $X$ is equal to  zero;
\item  a
 {\it generalized rotation}  if the weighted degree of $X$ is bigger than
 zero and less than one.
\end{enumerate}
\ed
 
 \bd(\cite {KMZ}){ We say that $X \in \aut(M_P,p)$  is an  $l$-integration
of a  rigid vector field $Y$  (necessarily in $\aut(M_P,p))$  if  the string of brackets  given by  $[ \dots [[X;\partial_{w}];\partial_{w}];\dots ];\partial_{w}]=Y$  is of length  $l.$ }   \ed

\section{Computing $\5g_c$}

Recall that  the component  $\5g_c$  of $\aut(M_H,0)$ is   by definition the set of (possible)  generalized rotations, that is, the set of (possible) rigid fields  of weight strictly bigger than $0.$ (See \cite{KM1}). As recalled in the introduction, this component is not trivial only in the singular Levi form case. We refer the reader to \cite{FM} for a "model"  example that illustrates   this phenomenon.
In this section we  derive an
 explicit description of  all model hypersurfaces  which  admit a nontrivial generalized rotation.
We start with the following lemmas.
 
\bl \label{whi}
Let  $a_j(z), j=1, \dots, n,$ be  $n$  $ \Bbb C-$ linearly independent holomorphic   polynomials  in the variable $z \in \Bbb C^N.$ Then 
there exist $z_1, \dots, z_n \in \Bbb C^N$ such that the determinant of the matrix 
$$
\begin{pmatrix}
a_1(z_1)&a_2(z_1)&\dots&a_n(z_1)\\
\vdots \\
a_1(z_n)&a_2(z_n)&\dots&a_n(z_n)\\
\end{pmatrix}
$$

is non zero.

\el
\begin{proof}
Let $V_z$ be the complex hyperplane in $\Bbb C^n$ given by  $$V_z= \{ (\lambda_1, \dots, \lambda_n) \in \Bbb C^n |\sum \lambda_j a_j(z) =0. \}$$ Consider the intersection $V$ of those varieties $V_z$ indexed by $z.$  Then, for each compact $Q$, we have 
$$ V\cap Q = (V_{z_1}\cap\ \dots\cap V_{z_k})\cap Q$$ for some $k.$ Since $V=\{0\},$   $$ (V_{z_1}\cap\ \dots\cap V_{z_k})\cap Q =\{0\}$$  for $Q$ containing $0.$ (See  Theorem 9C, page 100 in \cite {W}). Since  $V_z$ is a complex hyperplane, we  then obtain
 $$ V_{z_1}\cap\ \dots\cap V_{z_k} =\{0\}.$$ Then there need to exist $n$ of the $z_j$ with this property.   This achieves the proof of the lemma.
\end{proof}
\bl \label{do0} Let  $V_k, \ k\in \mathbb N$, be the  space
\begin{equation}
V_k=\{X |Y^k(X)=0                 \},
\end{equation}
where $X$ is a weighted homogeneous   holomorphic polynomial of a given   weighted degree  and $Y$ is a weighted homogeneous holomorphic vector field. Suppose that $V_1$ is not trivial. Then 
 there are strictly   positive  integers  $d_k \le k$ and $ g_k \le \dim V_k ,$ and a basis  of the form
\begin{equation}\label{do} \{ {F_s^k}\in   V_k, s=1, \dots, \dim V_k,|
\end{equation} $$
 \ \  Y^{d_k}({F_{s}^k}) =0,  \  Y^{{d_k}-1}({F_{s}^k})\ne 0, \  s= 1, , \dots, g_k, \  Y^{d_n-1}({F_{s}^n})=0,   s> g_k \} $$
such that   $\{Y^{{d_k}-1}({F_{s}^k})\}_ {s= 1}^ {g_k}$ are linearly independent.
\el
\bpf We prove the lemma by induction. Since $V_1$ s not trivial, the case $k=1$ is clear, with $d_1 =1.$ Suppose the lemma true for $k.$ We have $0=Y^{k+1}(X)= Y^k(Y(X))=0.$ The conclusion follows.
\epf
\bt 
\label{2piva} Let 
 $M_H$  be  given by  \eqref{2} admitting a generalized rotation
 $Y.$  
Then $P_C$ can be decomposed in the following way
\begin{equation}P_C=\sum_{j=1}^M T_j,\label{Barou19}
\end{equation} where each $T_j$ is given by 
\begin{equation}\label{Barou20}
 T_j =\Re (\sum_{k=1}^{N_j}<
  {U_j^{{(k)}}},
  { {{V_j^{{(N_j -k +1)}}}}}>), 
  \end{equation}
  where  $\{{ {{U_j^{(1)}}, \dots, {U_j^{(N_j)}} }}\}$  and  $\{{ {{V_j^{(1)}}, \dots, {V_j^{(N_j)}} }}\}$ are a symmetric pair of $Y_{s_j}-$ chains,  and $<\ ,\,>$ is the usual scalar product in $\Bbb C ^{s_j}.$
\et
\bpf
Let
\begin{equation}
P_C = \sum_{k=1}^l P_{c_k},
\end{equation}
where $P_{c_1} \neq 0, P_{c_l} \neq 0$, be the bihomogeneous expansion of $P_C$. Each $P_{c_j}$ is weighted homogeneous with respect 
to $z$ of weighted degree $c_j$  where $c_1< c_2 < \dots <c_l.$ 

We may write 
\begin{equation}
P_{c_1} =\sum_{j=1}^{r_1} {Q_j^{c_1}}{\overline {{Q_j^{\hat c_1}}}},
\end{equation}
with $r_1$ minimal. Since $Y$ is a generalized rotation, we must have 
\begin{equation}\label{do1} \overline Y(\sum_{j=1}^{r_1} {Q_j^{c_1}}{\overline {{Q_j^{\hat c_1}}})=
\sum_{j=1}^{r_1} {Q_j^{c_1}} \overline Y ({\overline {Q_j^{\hat c_1}}}} )= 0.
\end{equation}
Since $r_1$ is minimal, we have that $\{{Q_j^{ c_1}}\}_{j=1}^{r_1}$ are linearly independent and hence
\begin{equation}\label{do2}
 Y ( Q_j^{\hat c_1})= 0
 \end{equation}
 for all $j$.  
 We may assume that, 
  after a possible linear transformation     $$\{Y({Q_j^{c_1}}), \ Y({Q_j^{c_1}})\ne 0\} $$ are linearly independent. Let  $J_1= \{j |Y({Q_j^{c_1}})\ne 0\},$ and $J_2= \{j |Y({Q_j^{c_1}})= 0\}.$
  We may then rewrite $P_{c_1}$ as 
\begin{equation}
P_{c_1}=\sum_{j\in J_1}{Q_j^{c_1}}{\overline {{Q_j^{\hat c_1}}}} + \sum_{j\in J_2} {Q_j^{c_1}}{\overline {{Q_j^{\hat c_1}}}}
\end{equation}
We consider  the following subset of the set $\{P_{c_k}\},$ namely 
$P_k:=P_{c_1 +(k-1) \mu,}$ where $\mu>0$ is the weight of $Y.$
We claim that there exists $N\le l,$  such that $P_k, \ k \le N,$ can be  written as 
\begin{equation}\label{dodo}
P_{k} =\sum_{j=1}^{R_k}{Q_j^{c_k}}{\overline {{Q_j^{\hat c_k}}}}
+\tilde {P_{k}}
\end{equation}
so that
\begin{itemize}
\item $Y({Q_j^{c_{N}}})=0,$
\item $\{Y({Q_j^{c_k}})| \ Y({Q_j^{c_k}})\ne 0\} $ are linearly independent,
\item  $\{{Q_j^{ c_k}}\}_{j=1}^{r_k}$ are linearly independent,
\item   there is  $d_k$ such that 
$
  \overline{ {Y}^{d_k}}(\overline{Q_j^{\hat c_k}})=0,
$
 $\{\overline{ {Y}^{d_k-1}}
  (\overline{Q_j^{\hat c_k}})\ |\  \overline{ {Y}^{d_k-1}}
  (\overline{Q_j^{\hat c_k}})\ne 0\}$ are linearly independent, and $\overline{{Y}^{d_k-1}}(\tilde {P_{k} })=0.
$
\end{itemize} Note that $N$ is well defined  since $Y$ is a generalized rotation.
We prove the claim by induction. The case $k=1$ has just been proved.
Suppose  then  by induction that \eqref{dodo} holds for $k<N.$ 
We  write 
\begin{equation}
 P_{k+1} =\sum_{j=1}^{r_{k+1}} {S_j^{c_{k+1}}}{\overline{S_j^{\hat c_{k+1}}}}
\end{equation}
with $r_{k+1}$ minimal.

 Since $Y$ is a generalized rotation, we have
\begin{equation}\label{do5} \sum_{j=1}^{R_k}Y({{Q}_j^{c_k}}){\overline
 {{{Q}_j^{\hat c_k}}} + Y(\tilde {P_{k} }) +\sum_{j=1}^{r_{k+1}} {S_{j}^{c_{k+1}}}
 \overline Y({\overline {S_{j}^{\hat c_{k+1}}}}})= 0.
\end{equation}
Applying ${\overline Y}^{d_k}$ to \eqref{do5}, we get
\begin{equation}\label{do6}
\sum_{j=1}^{r_{k+1}} S_{j}^{c_{k+1}} \overline{ Y^{d_{k}+1}}
(\overline{{S_{j}^{\hat c_{k+1}}}})= 0.
\end{equation}
Since $r_{k+1}$ is minimal,
\begin{equation}
 \overline{ Y^{d_{k}+1}}
(\overline{{S_{j}^{\hat c_{k+1}}}})=0
\end{equation}
for all $j$.
Using \eqref{do},  we may then  rewrite $P_{k+1}$  in the form given by \eqref{dodo}, with $d_{k+1} \le d_k+1.$  The claim is then proved.
\\

We consider the following set $E$ given by
\begin{equation}\label{no}
E:=\{{{Q}_j^{c_k}}{\overline {{{Q}_j^{\hat c_k}}}}, \ \  \  j=1, \dots, R_k, \  k=1,
 \dots,  N \}.
\end{equation}

We claim that the following holds for every element of $E.$
\begin{enumerate}
\item $d_{k+1}= d_k +1,$\\
\item $Y({{Q}^{(c_k)}})= A_k {{Q}^{(c_{k+1})}} ,$\\
\item $ Y({{{Q^{(\hat c_{k+1})}}}})=B_{{
{k+1}}} {
{{{Q_k}^{(\hat c_k)}}}} +R_k,$
where 
${Y}^{d_k-1}({R_{k} })=0.$\\
\end{enumerate}

Suppose that this is true for $k < N-1$ 
and show that it is also true for $k+1.$ Using the fact that $Y$ is a generalized rotation, we have as in \eqref{do5}

\begin{equation}\label{do7} \sum Y({Q_j^{c_k}){\overline {{Q_j^{\hat c_k}}}} + Y(\tilde P_{k} ) +
\sum ({{Q}_{{j}}^{c_{k+1}}})\overline Y ({\overline {{Q}_{{j}}^{\hat c_{k+1}}}}}) + \overline Y(\tilde P_{k+1} )=0.
\end{equation}
Applying $ \overline Y^{d_{k-1}}$ to \eqref{do7}, we obtain, since $d_{k+1} \le d_k +1,$
\begin{equation}\label{do9} \sum Y({Q_j^{c_k}})  \overline Y^{d_{k-1}}({\overline {{Q_j^{\hat c_k}}}}) 
+ \sum ({{Q}_{{j}}^{c_{k+1}}})\overline Y^{d_k} ({\overline {{{Q}_{{j}}^{\hat c_{k+1}}}}}) =0.
\end{equation}
Hence, using \eqref{do9}, $d_{k+1}= d_k +1,$ and therefore, using \eqref{dodo} and  Lemma \ref{whi}
\begin{equation}
Y({Q^{(c_k)}})= A_{k} {Q^{(c_{k+1})}}.
\end{equation}
\begin{equation}
 Y^{d_k}({{{{Q}^{(\hat c_{k+1})}}}})=B_{{
 {k+1}}}Y^{d_{k}-1}{{{{Q}^{(\hat c_k)}}}}, 
\end{equation}
which implies 
\begin{equation}\label{do10}
 Y^{d_{k-1}}(Y({{{{Q}^{(\hat c_{k+1})}}}})-B_{{
 {k+1}}}{{{{Q}_k^{(\hat c_k)}}}})=0,
\end{equation}
and hence
\begin{equation}\label{do11}
 Y({{{{Q}^{(\hat c_{k+1})}}}})=B_{{
 {k+1}}} { {{{Q}^{(\hat c_k)}}}} +R_k,
\end{equation}
where 
$ {Y}^{d_k-1}({R_{k} })=0.$ This achieves the proof of the claim.
Using \eqref{do11} and \eqref{dodo}, we may then assume without loss of generality that $R_k=0.$
We define the  chains by putting 
\begin{equation}\label{no1}
\begin{cases}
 {U}^{(k)}_{1}  :={Q}^{(c_k)},\\
V_1^{{{(k)}}}  :={Q}^{(\hat c_{N-k+1})},\\

\end{cases}
\end{equation}
It follows from the above  properties of $E$ that $U_1^{(k)}$ and $V_1^{(k)}$
form a chain.   

In other words, we may write
\begin{equation}
P=\Re (\sum_{k=1}^{N}
 < {U_1^{(k)}},
  {{{V_1^{(N -k +1)}}}}>) + {\hat P}. \ \ \ 
\end{equation}

It follows from \eqref{do9} that $Y$ is a generalized rotation for  $$ \Im w =\Re (\sum_{k=1}^{N}
 < {U_1^{(k)}},
  { {{V_1^{(N -k +1)}}}}>).$$
  It follows from \eqref{do7} that $A_k = - ^t\overline { B_{k+1}}$, which means 
that the $U$ and $V$ are a pair of symmetric chains.
Hence $Y$ is a generalized rotation also for $\hat P$. We can repeat the above argument for 
$\hat P$ and in a finite number of steps we reach  the conclusion of the theorem.

\end{proof}
As  noticed in \cite{KM1},  symmetric chains and pairs of chains of any length can arise.

 \section{Computing  $\mathfrak{g_{nc}}$}

Recall that   the component  $\5g_{nc}$  of $\aut(M_H,0)$ is   by definition the set of (possible)  $1$-integration  of (possible) nontransversal shifts; they are   of weight strictly bigger than $0,$   defined up to $\5g_{c}.$ In  \cite{KM2}, we consider the case $n=2,$ and show that "only"  two model hypersurfaces occur for which $\5g_{nc}\ne 0,$ one being the "model" example studied in \cite{FM}    given by
$\Im w =  \Re z_1 \bar z_2^{d-1} .
$
In this section we  derive an
 explicit description of  all model hypersurfaces  which  admit a nontrivial $1$-integration and show that there exists a canonical representation of such a vector field.

Let 
$M_H$  of finite Catlin multitype  in $\mathbb C^{n+1}$ be  given by 
\begin{equation}\label{2}
M_H:=\{\Im w = P_C(z,\bar z)\}, \quad (z,w)\in \C^n\times\C,
\end{equation}
where $P_C$ is a weighted homogeneous polynomial of degree one with respect to the multitype weights $\mu_1, \mu_2, \dots, \mu_n.$ 
Suppose that   $M_H$ has  a (nontrivial)  $\mathfrak{g_{-\mu_l}},$ with $\mu_l$ chosen to be  minimal such that there exists a (nontrivial)  $X\in \mathfrak{g_{-\mu_l}}$ that admits a nontrivial $1$-integration.
By Lemma 6.1 in \cite{KMZ} there exist local holomorphic
coordinates preserving the multitype (with pluriharmonic terms
allowed), such that
\begin{equation}
  X=i \partial _ {z_{l}}.
  \label{xi}
\end{equation}
Hence we may write  ${P_C}$ in the following form

\begin{equation}
{P_C}(z, \bar z) = \sum_{j=0}^{m}{(\Re z_l)}^j P_j(z', \bar z'),\ \
 \label{pzz}
\end{equation}
for some homogeneous 
polynomials  ${P_j}$
in the variables $z' = ( z_1, \dots, \hat z_{l},  \dots, z_n)$,
with $P_m \neq 0.$
\bp
\label{first-nt}
Let $M_H$ be a holomorphically nondegenerate model with  $\5g_{nc}\ne 0.$ Let  $X$ be  given by
\eqref{xi}, 
and  {$P_C$  and $m$  be given by }\eqref{pzz}.
Then{
\begin{itemize}
\label{hopo}
\item  $ m \leq 2$
\item If $m=2,$  $P_2$ is a real constant.
\end{itemize}}
  
\ep
\begin{proof}
Splitting  $Y$ with respect to the powers of $z_l$, we obtain
\begin{equation}\label{yy0} Y= iw\partial_{z_l}+\sum_{j=-m}^k{}Y_j,
\end{equation}
where $Y_j$ is of the form
\begin{equation}\label{yy1}
Y_j={\phi_l}^j(z') {z_l}^{j+1}\partial_{z_l}+ < \phi^j(z')
{z_l}^{j}, \partial_{z'} >
+{\psi}^j(z')
{z_l}^{j+m}\partial_{w}.\end{equation}
(The coefficients are
zero when the power of $z_l$ is negative and  $Y_k \neq 0$).

We claim  that {$m-1\le k.$} Indeed, if not,  applying $\Re Y$ to
the first term  of the right handside of
{\eqref{pzz}},{ we obtain a term which contains } the
maximal{ nonzero  power }in $z_l$, namely
$$ - \frac{m}2 P_m^2{{(\Re z_l)}}^{2m-1}$$
while
all  other terms are of maximal power $m+k$
with respect to  $z_l$, {which gives  the  contradiction}.
Suppose by contradiction that $m>2.$
 We define $Z\in {\text{\rm aut}}(M_H, 0)$ by 
\begin{equation}\label{ouf}
Z:=[[Y,X],Y].
\end{equation}
$Z$ is nonzero non rigid vector field  since  $k\ge 2,$  and  its   weighted   homogeneous degree is  $2-3 \mu_l.$ By Theorem 1.3 in \cite{KMZ}, it implies that {$2-3 \mu_l \le1-\mu_l$}  by minimality of $\mu_l,$ since    $Z$  is not a $2-$integration of $\partial_{w}.$  Hence $\mu_l= \dfrac{1}{2},$ which is a contradiction with $m>2.$ This achieves the proof that $m \le 2.$

We now prove that $P_2$ is constant.
Let
$k$ be given by \eqref{yy0}. 
Let us denote the middle term of \eqref{yy1}
\begin{equation}\label{yy2}
 Y'_j :=< \phi^j(z') {z_l}^{j}, \partial_{z'} >
\end{equation}

and analogously
\begin{equation}\label{yy3}
 Y':= \sum_j Y'_j.
\end{equation}

Note that by weighted homogeneity
of $Y$, each coefficient  of {$Y'$} is weighted homogeneous in $z'$.

 First assume $k=1$.
From coefficients of degree three   with respect to $z_l,$
we obtain,
\begin{equation}\label{ho2aa}
 2x_l^{3} P_{2}^2 =  2x_l P_2  \Re \phi_l^1 z_l^2  + 2 x_l^2 \Re  Y_1'(P_2)
- \Im \psi^1 z_l^3.
\end{equation}

If $\psi^1(z') \ne 0$ in \eqref{ho2aa}, then it is a constant,
by comparing terms in $y_l^{3}.$
Hence, by homogeneity, $P_2$ is constant.
Next, assume  that  $\psi^1(z') = 0.$
Comparing degrees in $z'$, we see that $ \phi_l^1 $ and  $ P_2$ have the same
degree, or $\phi_l^1 =0. $ If $\phi_l^1 \neq 0, $then from the
coefficients of $x_l y_l^{2}$ we obtain that $\phi_l^1 $ is a
constant, hence $P_2$ is a constant.  On the other hand,
$ \phi_l^1 =0$ implies
\begin{equation}\label{ho2}
 x_l^{3} P_{2}^2 =  x_l^2 \Re Y_1'(P_2)
\end{equation}
which is impossible, unless $P_2 = 0$, since by positivity of $P_2^2$, the left hand side contains a non zero diagonal term in $z'$, while
the right hand side has no such terms.

Now assume that $k \ge 2$ (note that $k=0$ is impossible, since  $k \geq
m-1$). $Z$ given by \eqref{ouf}  is then a nonzero nonrigid  vector field, and then $\mu_l =\dfrac{1}{2},$ which means that $P_2$ is constant. This achieves the proof of the proposition. 
\end{proof}
We have the following lemma.
\bl
\label{lemtub}  Let $m \in \Bbb N, \ m {\ge 1}.$ There
exist uniquely determined nonzero complex numbers ${\alpha_0}, \dots,
\alpha_{m-1}$ such that for every $ z \in  \Bbb C,$
\begin{equation}
{{(\Re z)}^{2m-1} = \sum_{j=0}^{m-1} {(\Re z)}^{j} \Re (\alpha_j z^{2m-1-j}).}
\end{equation}

\el

\begin{proof} Indeed, by comparing coefficients of {$ z^{m-1}\bar z^{m} $}
we obtain the value of $\alpha_{m-1}$. Continuing this way, from the
coefficients of {$ z^{m-1-j} \bar z^{m+j} $ }we obtain the uniquely
determined values of $\alpha_{m-1-j}$.
\end{proof}
\bp
\label{fourth-nt}
Let $M_H$ be a holomorphically nondegenerate model,
with $P_C$  given by \eqref{pzz},  and  $m =2$.
Let $X=i \d_{z_{l}}$ be in $aut(M_H, 0).$ 
Then
there is a vector field $Y$ in $aut(M_H, 0)$ such that
$[Y,W]=X$, if and only if
 $P_C$ is biholomorphically equivalent, by a change of multitype coordinates, to
\begin{equation}
 P_C(z, \bar z) =
 x_l^2 +
 P_0( z', \bar z'),
\end{equation}
where  $P_0( z', \bar z')$ is a balanced polynomial without pluriharmonic terms.

Moreover $Y$ can be chosen canonically as
\begin{equation}
Y= iw\partial_{z_l}+ a{z_l}^2\partial_{z_l}  + z_l S  +  b{z_l}^3\partial_{w},
\end{equation}
where  $a$ and $b$ are uniquely determined nonzero constants, 
$S= < \phi(z')
, \partial_{z'} >$
 uniquely determined by the condition 
 $S(P_0)=P_0.$
\ep
\begin{proof} Let $Y$ be given by \eqref{yy0} and \eqref{yy1}.
Without  loss of generality, we may assume that both $P_1$ and $P_0$
contain no pluriharmonic terms.
Indeed  pluriharmonic terms in $P_1$ can be eliminated
by a change of variables
$z_l^* = z_l + S(z')$, where $S$ is a holomorphic polynomial in $z',$ using the fact that $P_2$ is constant. Then to eliminate pluriharmonic  terms in $P_0,$ we perform a change of coordinates of the form $w^* = w + H(z'),$ where $H$ is a holomorphic polynomial in $z'.$
 We claim that  $P_1 = 0$.
Applying $\Re Y$ to $P_C-v$ gives
\begin{equation}\label{hahu}
 -(2 x_l +  P_1 )( x_l^2 + x_l P_1
 + P_0)
 + 2 \Re( \phi_l \frac{\partial P_C}{\partial z_l}) + 2 \Re
(\sum_{j\neq l} \phi_j \frac{\partial P_C}{\partial z_j} )-  \Im \psi = 0.
\end{equation}
Let $k$ be as in \eqref{yy0}.
 Assume first that
$k=1$. For the third order terms in $z_l$ we obtain
\begin{equation}
2 x_l^3  = 2x_l \Re ( \phi_l^{1} z_l^{2})
 - \Im (\psi^{1}z_l^{3}).
\end{equation}
By Lemma \eqref{lemtub}, $\phi_l^{1}$
and $\psi^{1}$  are  unique  non zero constants, $\phi_l^{1}\in \Bbb R^*.$
Looking at terms of second order in $z_l$ we obtain from \eqref{pzz}
and \eqref{hahu}
\begin{equation}
 - 3x_l^2 P_1
 + x_l  \Re ( \phi_1^0 z_l)
 + 2x_l \Re( \sum_{j\neq l}^n \phi_j^1 z_l \frac{\partial P_{1}}{\partial z_j}) + \Re (\phi^1_l z_l^2 P_1) -
 \Im \psi^{0}z_l^{2} =0.
\end{equation}
Looking at coefficients of $y_1^2$, we obtain that $P_1$ is pluriharmonic, since $\phi_l^{1}\in \Bbb R^*$.  Hence $P_1 =0$.
{
Observe that this implies  
\begin{equation}\label{id1}
Y_0 \in \5g_c.
\end{equation}
Next, let  $k>1$. Then  $k=2$ since $\mu_l=\dfrac{1}{2}.$  We may assume that after a linear change of coordinates,
\begin{equation}
Y_2 ={Y_2}'=  z_l^{2}\frac{\partial }{\partial z_s}, \ s\ne l.
\end{equation}
We obtain, as before, for the third order terms in $z_l$ 
\begin{equation}
2 x_l^3  =2 x_l \Re ( \phi_l^{1} z_l^{2}) + 2x_l \Re(  {z_l}^2 \frac{\partial P_{1}}{\partial z_s})
 - \Im (\psi^{1}z_l^{3}).
\end{equation}
By Lemma \eqref{lemtub} and the fact that $P_1$ does not contain harmonic terms, we obtain  $\frac{\partial P_{1}}{\partial z_s}=0,$ and hence  $\phi_l^{1}$
is a  non zero real   constant. In fact one may compute directly that $\phi_l^{1}=\dfrac{3}{2}.$
For the terms of second order in $z_l,$ we get
\begin{equation}\label{chaban1}
 - 3x_l^2 P_1
 + x_l  \Re ( \phi_1^0 z_l)
 + 2x_l \Re( \sum_{j\neq l}^n \phi_j^1 z_l \frac{\partial P_{1}}{\partial z_j})\end{equation} $$+ \  \Re (\phi^1_l z_l^2 P_1) + 2\Re( {z_l}^2 \frac{\partial P_{0}}{\partial z_s}) -
 \Im \psi^{0}z_l^{2} =0.
$$
Looking at coefficients of $y_1^2$, we obtain that 
$P_0$ contains a non zero term of the form $z_s H,$ where $H$ does not depend on the variable $z_s$ and is of the same degree as $P_1$ with
\begin{equation} \label{fer3}
\dfrac{3}{2} P_1 + H +\bar {H}  =0.
\end{equation}
 By  definition of the weights, since $z_s$ has weight 
$\mu_s=\frac{1}{2},$ that forces $H$ to be linear, and hence, using \eqref{fer3},  $P_1$  should contain the harmonic  term $H,$ which is impossible. It shows that $k=2$ is impossible.
}

Returning to the only possible case, k=1, and
looking at the linear terms in $z_l,$ we obtain
\begin{equation}
- 2 x_l P_0  + 2 \Re ( z_l  \sum_{j\neq l}^n
\phi_j^1 \frac{\partial P_0}{\partial z_j}) + x_l \Re (\phi_1^{-1}) -  \Im (\psi^{-1} z_l)= 0
\end{equation}
which gives equations for coefficients of $x_l$ and $y_l$.  Namely
\begin{equation}
 -2 P_0  + 2 \Re ( \sum_{j}^n
\phi_j^1 \frac{\partial P_0}{\partial z_j}) +   \Re (\phi_1^{-1}) -  \Im (\psi^{-1})= 0
\end{equation}
and
\begin{equation}
   - 2\Im  \sum_{j=2}^n
\phi_j^1 \frac{\partial P_0}{\partial z_j} -  \Re \psi^{-1}= 0.
\end{equation}
Using
 the fact that $P_0$ contains no pluriharmonic terms, and that the weight of  $\psi^{-1}$ is equal to one,  it follows that  $\psi^{-1}=0 $.
 By the same argument, $\phi_1^{-1}=0$.  Hence we obtain   
\begin{equation}
   \Im  \sum_{j=2}^n
\phi_j^1 \frac{\partial P_0}{\partial z'} = 0,
\end{equation}
hence
\begin{equation}
 P_0 =  \sum_{j=2}^n
 \phi_j^1 \frac{\partial P_0}{\partial z'}.
\end{equation}
It follows that $P_0$ has a complex reproducing field, hence $P_0$
is a balanced polynomial,
as claimed.
That finishes the proof.
\end{proof} 

Now we consider the case $m=1$.
Let us write 
$P_0$ as
\begin{equation}\label{mmc1}
 P_0 = \sum_{j=2}^{s} S_j(z') \overline{ Q_j(z')},
\end{equation}
where  $s$ is minimal. Without loss of generality, we may assume that  $P_0$ contains no pluriharmonic terms
  by performing local holomorphic change of coordinates.

\bd
For a (n-1)-tuple of functions $R= (R_1, \dots R_{n-1})$ depending on $n-1$ complex variables $Z_1, \dots ,Z_{n-1},$
we  denote by $\Delta (R_1, \dots R_{n-1})$ the determinant of the Jacobi matrix of $R_1, \dots,  R_{n-1}$
$$ \Delta(R_1, \dots R_{n-1}) 
= \det ( \frac{\partial R}{\partial {Z}_1}, \dots   \frac{\partial R}{\partial {Z}_{n-1}})$$

and set

$$ {\Delta_j}^H(R_1, \dots R_{n-1})  = \det ( \frac{\partial R }{\partial {Z}_1}, \dots,  \frac{\partial R  }{\partial {Z}_{j-1}},  H,  \frac{\partial R  }{\partial {Z}_{j+1}},  \dots, \frac{\partial R}{\partial {Z}_{n-1}}).$$
\ed

\bp
\label{third-nt}
Let $X=i \d_{z_{l}}$
be in
$ aut(M_H, 0)$ and $P_C$ a homogeneous polynomial of degree $d$ of the form \eqref{pzz} with $m =1$.  Assume that $P_0$ is given by \eqref{mmc1}.
Then there  is a vector field $Y$ in $aut(M_H, 0)$ such that
$[Y,W]=X$ if and only if
\begin{enumerate}
\item $ P_1 = \Re (Q_1)$
where $Q_1$ is a holomorphic polynomial, 
\item  for every choice of  $ Q:=(Q_{j_1}, \dots, Q_{j_{n-1}})$ such that $\Delta (Q_{j_1}, \dots, Q_{j_{n-1}}) \ne 0,$ there exist  homogenous functions $g_k, \ \ k=1, \dots, s,$ of degree one, holomorphic outside a  analytic set, such that  $Q_k= g_k(Q_{j_1}, \dots, Q_{j_{n-1}}),\ \   k=1, \dots, s$ 
\item the polynomials  $$\frac{1}{2}Q_1 {\Delta_j}^{Q}(Q_{j_1}, \dots, Q_{j_{n-1}})$$   are divisible by  $ \Delta(Q_{j_1}, \dots, Q_{j_{n-1}}) $ 
in the ring of  holomorphic polynomials.
\end{enumerate}
{Moreover $Y$ can be chosen canonically as
\begin{equation}
Y= iw\partial_{z_l}+ {\phi}_0 {z_l}\partial_{z_l} + \psi {z_l}\partial_{w} + S
\end{equation}
where  ${\phi}_0$ and $\psi $  are  homogeneous polynomials  in the variables $ z'$  uniquely determined  and satisfying   the condition   $\ i{\phi}_0 \partial_{z_l} + i\psi \partial_{w}\in \5g_c ,$   
$S= < \phi(z')
, \partial_{z'} >$
 uniquely determined by the condition 
 $2S(P)= Q_1P.$}
\ep

Let us remark that in complex dimension three, the previous proposition implies that 
$P_0 =0$, which was already proved in \cite{KM2}. 

\begin{proof}
Integrating $X$, we obtain the same form of $Y$ as before, 
\begin{equation}
 Y = i  w \d_{z_l} + \sum_{j=1}^n \phi_j \d_{z_j} + \psi \d_w.
\end{equation}
From $\Re Y(P-v) = 0$, using $\Re X(P) = 0$, we obtain

\begin{equation}
 P_0 P_1 + x_l P_1^2 = 2 x_l \Re \sum_{j \neq l} \phi_j \frac{\partial P_{1}}{\partial z_j} +
  \Re \phi_l
 P_1 + 2 \Re \sum_{j \neq l} \phi_j \frac{\partial P_{0}}{\partial z_j}
-  \Im \psi.
\end{equation}
Let $k$ be again as in \eqref{yy0}. 
For the constant and linear terms in $z_l$ we have
\begin{equation}
 P_0 P_1 = 2 \Re   \sum_{j \neq l} {\phi_j}^0 \frac{\partial P_{0}}{\partial z_j}
+  \Re \phi_l^{-1} P_1 -
\Im \psi^{-1}
\label{x10}
\end{equation}
and
\begin{equation}
 x_l P_1^2 = 2 x_l \Re \sum_{j \neq l} {\phi_j}^0 \frac{\partial P_{1}}{\partial z_j}
 + \Re \phi_l^0 z_l  P_1 + 2 \Re z_l \sum_{j \neq l} {\phi_j}^1 \frac{\partial P_{0}}{\partial z_j}-
\Im \psi^0 z_l.
\label{x11}
\end{equation}
Let
first
$k =2$. By \eqref{ouf}, $\mu_l=\dfrac{1}{2}.$ 
and hence $P_1= \Re H,$ where $H$ is  holomorphic linear, { using the definition of weights.} 
 After performing  local holomorphic changes of coordinates, we may assume that  $H=z_k$  for some $k\ne l.$  { Because of minimality in \eqref{mmc1}, we can normalize $P_0$ and assume that $H=z_k$ is not in the $\Bbb C$ linear span of  the $S_j$ 
by absorbing such a term into $P_1.$  }
We get

\begin{equation}
 0 = \Re (\phi_l^{2} z_l^{3} (\frac{z_k +\overline {z_k}}{2})) +  x_l \Re ({z_l}^2
 {\phi_k}^2)
 -  \Im (\psi^{2}z_l^{3}).
\end{equation}
From  the coefficient of $\bar z_l z_l^{2}$, we obtain  that ${\phi_k}^2=0.$
If  $\phi_l^{2}\ne 0,$ it follows that
 $\Re \phi_l^{2} z_l^{3} P_1$ is pluriharmonic,
hence
$P_1$
is constant, which is impossible.
We  then assume $\phi_l^{2}= 0.$ Looking for terms of second degree in $z_l,$ we obtain
\begin{equation}\label{ret}
0 = \Re (\phi_l^{1}z_l^{2} (\frac{z_k +\overline {z_k}}{2})) + x_l\Re (z_l  {\phi_k}^1 )) - \Im (\psi^{1}z_l^{2}) + 2\Re ({z_l}^2
 \sum_{j \neq l} {\phi_j}^2 \frac{\partial P_{0}}{\partial z_j}).
\end{equation}
Looking at the term $z_l\bar z_l$ in \eqref{ret}, we obtain that ${\phi_k}^1=0,$  since ${\phi_k}^1$ has weight $\dfrac{1}{2}.$
 If $\phi_l^{1}\ne 0,$ then we get a contradiction since {$z_k$ is not in the $\Bbb C$ linear span of  the $S_j.$ }
 If $\phi_l^{1}= 0,$  then 
\begin{equation}
\Re ({z_l}^2\sum_{j \neq l} {\phi_j}^2 \frac{\partial P_{0}}{\partial z_j})=0,
\end{equation}
and hence 
\begin{equation}
\sum_{j \neq l} {\phi_j}^2 \frac{\partial P_{0}}{\partial z_j}=0,
\end{equation}
Using the fact that $M_H$ is holomorphically nondegenerate, we obtain that $$\sum_{j \neq l} {\phi_j}^2 \frac{\partial }{\partial z_j}=0,$$  since ${\phi_k}^2=0.$  This means that $k=2$ is impossible.

Now let $k=1$. { First of all, observe that $Y_1(P_1)=0.$}

Since $[Y,X]$ is a generalized rotation  ($d>2$), $P$ has a chain structure.

By the results of Section 3, we can write (in the scalar product notation),

$$ P = \sum Re < U_j, \bar U_{n-j+1} >,$$
where 
$$[Y,X](U_n) = 0.$$

Since the degree of $[Y,X]$ is $d-2$, 
we obtain that the {maximal} length of {a} chain is  two, hence the first element in those maximal  chains is linear, while the second one has degree $d-1$.
 { That means that  nonlinear   terms  could exist but are killed by   $[Y,X].$}
{We also notice that $[Y,X]= -i\phi_l^0\d_{z_l} -i \psi^0\d_{w}-iY_1,$  which means that $[Y,X](P_0)= -iY_1(P_0).$}
Write
$$ P_1(z', \bar{z'}) = \Re H(z') +\tilde {P}(z', \bar{z'}) , $$  where $H$ is harmonic.

 We  can write 
{
\begin{equation}
  P_{0}(z', \bar{z'}) = \Re (\sum_{j=1}^r L_j(z') \bar S_j(\bar{z'})) + \tilde{P_0}(z', \bar{z'}) ,
\end{equation}
where $Y(\tilde{P_0}(z', \bar{z'}))=0,$} and 
where $L_j$ are linear functions and $r$ is minimal. 
Because of minimality, we can normalize $P_0$ and assume that $H$ is not in the $\Bbb C$ linear span of  the $S_j$, 
by absorbing such a term into $P_1.$

Let  $$ Y_1 = \sum_{j \neq l} f_j \d_{z_j}.$$
From the coefficients of $y_l$ in \eqref{x11} we obtain
\begin{equation}
  \Im \phi_l^0   P_1 + 2 \Im Y_1 (P_{0})
+ \Re \psi^0 = 0. 
\label{x1122}
\end{equation}
{ We first assume that  $\tilde {P}(z', \bar{z'})=0,$ that is, $P_1$ is harmonic. 
Consider the diagonal terms in this equation and in the real part of \eqref{x11}.
From the real part, on the l.h.s we obtain just $H \bar H$. $P_0$ cannot produce such a term, since $H$ is not in the $\Bbb C$ linear span of  the $S_j$, we get 
 $ \phi_l^0 = cH$, $c \in \mathbb R$. Now looking at the imaginary part, the  first term compensates with the third one,  hence 
 $$\Im Y_1 (P_0)= 0.$$
 Since there are no more diagonal terms on the r.h.s of the real part of \eqref{x11}, we obtain  
 $$\Re Y_1 (P_0)= 0.$$
  It implies $ Y_1(P_0) = 0$, hence $M_P$ is holomorphically degenerate, which is a contradiction. 
  Assume now that assume that  $\tilde {P}(z', \bar{z'})\ne 0.$ Then \eqref{x1122} holds if $\phi_l^0=0.$ Using \eqref{x11}, we obtain that $H=0.$  Applying $Y_1$ to  \eqref{x10}, and using the fact that $Y_1(P_1)=0$ we obtain a contradiction.}
  It follows that $k=1$ is impossible.

Now, let $k=0$.
 We get
\begin{equation}
 x_l P_1^2 = 2 x_l \Re
  (\sum_{j \neq l} {\phi_j}^0 \frac{\partial P_{1}}{\partial z_j})
+
 \Re ( \phi_l^0 z_l P_1) - \Im ( \psi^{0} z_l).
\label{x11bb}
\end{equation}
From the coefficients of $y_l$, we get
\begin{equation} - \Im \phi_l^0 P_1 - \Re  \psi^{0}= 0. \label{yy11}
\end{equation}

 This implies  that $P_1$ is pluriharmonic, namely
$P_1 = c \Re \phi_l^0.$ 
Notice that
$ \phi_l^0 =0$ leads to contradiction. Indeed,
if $ \phi_l^0 =0$, then $\psi^0 =0$,
since $P_1$ cannot be constant.
It follows that
\begin{equation}\label{ho2}
  P_{1}^2 =  2 \Re
 \sum_{j \neq l} \phi_j \frac{\partial P_{1}}{\partial z_j}
\end{equation}
which is impossible, since the left hand side contains a nonzero balanced term
 in $z'$,
while
the right hand side has no such terms.
That gives the contradiction.

%
%
%

Next consider the equation for the coefficients of $x_l$ in
\eqref{x11bb},
\begin{equation}
  P_1^2 = 2  \Re Y'_0(P_1) +
 \Re \phi_l^0 P_1 - \Im  \psi^{0}.
\label{x111bbb}
\end{equation}
Substituting $P_1 = c \Re \phi_l^0$, from the mixed terms we obtain  $c=1$.

For terms of order zero we
obtain
\begin{equation}
  P_0 P_1=   \Re \phi_l^{-1} P_1 +
2 \Re Y'_0 (P_0)
- \Im \psi^{-1}.
 \label{newx11}
 \end{equation}
%
%
%

Let us write now
$P_0$ as
\begin{equation}
 P_0 = \sum_{j=2}^{s} S_j(z) \bar Q_j(z),
\end{equation}
where  $s$ is minimal as above.
From now on, we denote $Q_1={\phi_l}^0.$ 
From equations \eqref{x111bbb}, \eqref{newx11} we get the following equations
\begin{equation}\label{miaou1}
{Y_0}'(Q_1)= \dfrac{1}{2}{Q_1}^2,
\end{equation}
and
\begin{equation}\label{miaou2}
{Y_0}'(Q_j)= \dfrac{1}{2}{Q_1}Q_j.
\end{equation}

By holomorphic nondegeneracy and reality, the polynomials $Q_1,  Q_2, \dots, Q_s$
are generating, i.e. their gradients span $\mathbb C^{n-1}$ at a generic point. This gives $s \geq n-1$. We may then choose $ Q:=(Q_{j_1}, \dots, Q_{j_{n-1}})$ such that $\Delta (Q_{j_1}, \dots, Q_{j_{n-1}}) \ne 0,$
For every $k =1, \dots, s,$  there exist holomorphic  functions $g_k$ of $n-1$ variables in a neighborhood of a generic point,  
such that $Q_k= g_k (Q_{j_1}, \dots, Q_{j_{n-1}}).$ Substituting  for $Q_k$ into \eqref{miaou1}, \eqref{miaou2}, we obtain

$$ {Y_0}'(Q_k) = \nabla g_k ({Y_0}'(Q_{j_1}), \dots, {Y_0}'(Q_{j_{n-1}})) = 
\dfrac12  \nabla g_k (Q_1 Q_{j_1}, \dots,  Q_1 Q_{j_{n-1}}). 
$$ 

On the other hand, 
\begin{equation}
{Y_0}'(Q_k)= \dfrac{1}{2}{Q_1Q_k}.
\end{equation}
Hence
\begin{equation}
 g_k(Q) = <\nabla g_k,  Q>
\end{equation}
It follows  that $g_k$ is homogeneous of degree one. 
Now, in order to determine   the component
of $Y'_0$, we use   Cramer's rule. This leads to
$$ {\phi^0}_j =\frac{1}{2}Q_1\dfrac{ {\Delta_j}^{Q}(Q_{j_1}, \dots, Q_{j_{n-1}})}{ \Delta(Q_{j_1}, \dots, Q_{j_{n-1}})} $$  

This implies the statement of the proposition.
\end{proof}

We will now show that there is no analogous statement in the case of unequal weights. 

\be{\label{ex5}}
Consider a model in $\mathbb C^4 $, given by  
\begin{equation}
 P(z, \bar z) = x_1 Re i z_2^l  + S(z_3, \bar z_3) Re z_2^l,
\end{equation}
where $S$ is a homogeneous real valued polynomial in $z_3$ of degree bigger than one and $l$ is an integer.
Note that $P_0$ is not balanced for 
suitable $S$. More concretely, let us take $l = 3$ and 

\begin{equation}
 S(z_3, \bar z_3) = \Re z_3 \bar z_3^3.
\end{equation}
The multitype weights become $(\frac14, \frac14, \frac 1{16})$.

Taking
{\begin{equation}
 Y' = \frac{i}6 z_2^4 \partial_{z_2}, 
\end{equation}}
we obtain a symmetry which is not of the form described in the previous result.

\ee

\section{Proofs of the main results}

In this section we complete the proofs of the results stated in the introduction.

%
{\it Proof of Theorem \ref{th12}.}
We apply Propositions  \ref{first-nt},  \ref{fourth-nt} and  \ref{third-nt}.

{\it Proof of Theorem \ref{pivo}.} This is Theorem \ref{2piva}.

{\it Proof of Theorem \ref{emb}.}
Let $Y$ be a generalized rotation.
In the notation  of Theorem \ref{2piva}, we set
\begin{equation}
 K =  2 \sum_{j=1}^{M}s_j N_j+1.
\end{equation}
We define a hyperquadric in  $\mathbb C^{K+1}$ by
\begin{equation}
 \Im \eta = \Re \sum_{j=1}^{M}\sum_{k=1}^{N_j}
  < \zeta_{j,(k)},
  {\zeta'_{j, (N_j -k +1})}>, 
\end{equation}
and consider the mappping  $\mathbb C^{n+1} \to \mathbb C^{K+1}$
given by $\eta = w$ and 
\begin{equation}
 \zeta_{j,(k)} = U_j^{(k)}(z).
\end{equation}
and
\begin{equation}
 \zeta'_{j,(k)} = V_j^{(k)}(z).
\end{equation}
It is immediate to verify that the automorphism $Y$ of $M_P$ is $f$-related to the automorphism of this hyperquadric, defined by 
\begin{equation}
 Z =  \sum_{j=1}^{M}\sum_{k=2}^{N_j}
 <A_{{k-1}, j} \;  \zeta_{j,{(k)}} ,\overline{\d_{ \zeta_{j,(k-1) }}}> + <^tB_{{k-1}, j} \;  {\zeta}'_{j,{(k)}} ,\overline{\d_{ {\zeta}'_{j,(k-1) }}}>.
\end{equation}
Indeed, the condition for $f$-related vector fields becomes exactly the chain condition 
\eqref{yuv}-\eqref{dcj}.

If $\5g_1 \neq 0$, then by Theorem \ref{cap}
\begin{equation}
 P(z, \bar z) = \sum_{|\al|_{\Lambda}= |\bar \al|_{\Lambda} = 1} A_{\al, \bar \al} z^{\al} \bar z^{\bar \al},
\end{equation}
where $A_{\al, \bar \al}\ne 0.$
 We order the multiindices and write
$P$ as
\begin{equation}
 P(z, \bar z) =\Re( \sum_{j=1}^{R}  z^{\al_j}(\sum_{k=1}^{N_j} A_{j,k}
  \bar z^{\al_{k,j}})).
\end{equation}
Consider the hyperquadric in  $\mathbb C^{R + (\sum_{j=1}^R N_j)+1}$ defined by
\begin{equation}
 \Im \eta = \Re (\sum_{j=1}^{R}
   \zeta_{j}(\sum_{k=1}^{N_j} A_{j,k}
  \bar \zeta_{{k,j}}))
  ,
\end{equation}
and the mappping  $f : \mathbb C^3 \to \mathbb C^{R + (\sum_{j=1}^R N_j)+1}$
given by $\eta = w$ and
 $\zeta_{j} = z^{\al_j}$ for $ j = 1, \dots, R, \ \ \zeta_{{k,j}}= z^{\al_{k,j}}, \ k=1, \dots, N_j. $

It is immediate to verify that the vector field in $aut(M_H, 0)$
$$Y =\left(\sum_{j=1}^n \lambda_j z_j  \d_{z_j} \right) w+
w^2 \d_{w},\
$$
is $f$-related to the infinitesimal automorphism of the above hyperquadric
 given by
$$Z  = \eta (\sum_{j=1}^{R} \zeta_j  \d_{\zeta_j} +  \sum_{j=1}^{R}(\sum_{k=1}^{N_j} \zeta_{k,j } \d_{\zeta_{k,j}} )      )
+ \eta^2 \d_{ \eta}.\ $$

The case  $\5g_{nc}\neq 0$ is completely analogous.


\begin{thebibliography}{BER96b}
%
\itemsep=2pt

\bibitem{BEH}
Baouendi, M. S.
 Ebenfelt, P.,
Huang, X.
 Super-rigidity for CR embeddings of real hypersurfaces into hyperquadrics
Adv. Math.
{\bf 219}
(2008),
1427--1445

\bibitem{BP}
{Bedford, E.,  Pinchuk, S. I.,
 {\it Convex domains with noncompact groups of automorphisms},
  {Mat. Sb.},
    \textbf{185}
     (1994), 3--26.

\bibitem{BES}
Beloshapka, V. K., Ezhov, V. V.,  Schmalz, G.,
{\it Holomorphic classification of four-dimensional surfaces in
              {$\mathbb C^3$}},
 Izv. Ross. Akad. Nauk Ser. Mat.},
 \textbf{72} (2008), 3--18.


 \bibitem{Beloshapka}
Beloshapka, V. K., {\it Symmetries of real hypersurfaces in
complex 3-space}, Math. Notes, \textbf{78} (2005) 156--163.

\bibitem{BG}
Bloom, T.,  Graham, I., {\it On "type"  conditions for
generic real submanifolds of $  C\sp{n}$},  Invent. Math.  \textbf{40}
(1977),  217--243.

\bibitem{BFG} Beals, M., Fefferman, C., Graham R.,
{\it Strictly pseudoconvex domains in $\mathbb C^n$},
 Bull.\ Amer.\ Math.\ Soc.\ (N.S.)
 \textbf{8} (1983), 125--322.



\bibitem{C}
  Catlin, D.,
   {\it Boundary invariants of pseudoconvex domains},
Ann.\ Math.\ {\bf 120} (1984), 529--586.

\bibitem{C1}
  Catlin, D.,
   {\it Subelliptic estimates for $\bar \partial$-Neumann problem  on pseudoconvex domains },
Ann.\ Math.\ {\bf 126} (1987), 131--191.


\bibitem{CM} Chern, S.\ S.\, Moser, J., \textit{Real hypersurfaces in
complex manifolds},  Acta Math.\  \textbf{133} (1974),  219--271.



\bibitem{EHZ}
Ebenfelt, P.,
Huang, X.,
Zaitsev, D.,
The equivalence problem and rigidity for hypersurfaces embedded into
   hyperquadrics
 Amer. J. Math.
{\bf 127}
(2005)
169--191.


\bibitem{ELZ2} Ebenfelt, P., Lamel, B., Zaitsev, D., {\it
Finite jet determination of local analytic CR automorphisms and
   their parametrization by 2-jets in the finite type case},
  {Geom. Funct. Anal.},
   {\bf 13},
   (2003),
   3, 546--573.



\bibitem{FK}
Fels, G., Kaup, W., {\it Classification of Levi degenerate homogeneous
CR manifolds in dimension 5}, Acta Math. \textbf{201} (2008),  1--82.

\bibitem{HY} Huang, X.,  and Yin, W.,
{\it A Bishop surface with a vanishing Bishop invariant},
Invent. Math. {\bf 176} (2009), 461--520.



\bibitem{K} Kohn, J.\ J., \textit{Boundary behaviour of
$\bar \partial$ on weakly pseudoconvex manifolds of dimension
two},
 J.\ Differential  Geom.\  \textbf{6} (1972),  523--542.

%

\bibitem{KoLa} Kol\'a\v r, M., Lamel, B. {\it Ruled hypersurfaces in
$\mathbb C^{2}$}, Journal of Geometric Analysis, {\bf 25} (2015),
1240--1281.

\bibitem{KMZ} Kol\'a\v r, M., Meylan, F., Zaitsev, D.,  \textit{Chern-Moser operators and
polynomial models in CR geometry}, Advances in Math. {\bf 263}
(2014), 321--356.


\bibitem{KM1} Kol\'a\v r, M., Meylan, F., \textit{Higher order symmetries of real hypersurfaces in $\mathbb C^3$},
 Proc. Amer. Math. Soc. {\bf 144}
(2016), 4807--4818.
\bibitem{KM2} Kol\'a\v r, M., Meylan, F., \textit{Nonlinear CR automorphisms of Levi degenerate
hypersurfaces  and a new gap phenomenon},
 Annali della Scuola Normale Superiore di Pisa Cl. Sci. {\bf 19} (2019), 847-868 
\bibitem{AIM} AIM list of problems in CR geometry, available at \newline http://www.aimath.org/WWN/crmappings/crmappings.pdf



\bibitem{Ko1} Kol\'a\v r, M., {\em The Catlin multitype and biholomorphic equivalence
of models}, Int. Math. Res. Not. IMRN {\bf 18} (2010), 3530--3548.

\bibitem{Ko1a} Kol\'a\v r, M., {\it Normal forms for hypersurfaces of finite type in $ \mathbb C^2$}, 
Math. Res. Lett., \textbf{12} (2005),  897--910.

\bibitem{FM}
Kol\'a\v r, M., Meylan, F., {\it Chern-Moser operators and weighted jet determination problems},
Geometric analysis of several complex variables and related topics, 75--88, Contemp. Math. 550, 2011.

\bibitem{KS}
Kossovskiy, I., Shafikov, R., {\it Analytic differential equations
and spherical real hypersurfaces}, Journal Diff. Geom., \textbf{102}
(2016), 67--126.

\bibitem{KZ}
Kossovskiy, I., Zaitsev, D., {\it Convergent normal form for real hypersurfaces at generic Levi degeneracy}, arXiv:1405.1743
(2014).

\bibitem{KL}
Kruzhilin, N. G., Loboda, A. V.,
 {\it Linearization of local automorphisms of pseudoconvex surfaces},
 Dokl. Akad. Nauk SSSR,
\textbf{271} (1983), 280--282.


\bibitem{Po} Poincar\'e, H.,  \textit{Les fonctions analytiques de
deux variables et la repr\'esentation conforme, } Rend. Circ. Mat.
Palermo \textbf{23} (1907), 185--220.


\bibitem{S} Stanton, N.,
{\it Infinitesimal CR automorphisms of real hypersurfaces}, Amer. J. Math. {\bf 118} (1996),  209--233.


\bibitem{V} Vitushkin, A.G., : \textit{Real analytic
hypersurfaces in complex manifolds}, Russ. Math. Surv. \textbf{40}
(1985),  1--35.
\bibitem{W} Whitney, H., : \textit{Complex analytic varieties}, Addison-Wesley Series in Mathematics. 
(1972).






%
%
%
%
%
%
%
%
%
%
%
%
%
%
%
%
%
%
%
%
 \end{thebibliography}
\end{document}